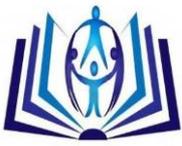

# Regional Boundary Strategic Sensors Characterizations

Raheam A. Al-Saphory[1], Hind K. Kolaib[2]

[1,2]Department of Mathematics, College of Education for Pure Sciences, Tikrit University, Tikrit, IRAQ.

[1]E-mail: saphory@hotmail.com and, [2]Hind.K.Kolaib@hotmail.com

**ABSTRACT**

This paper, deals with the linear infinite dimensional distributed parameter systems in a Hilbert space where the dynamics of the system is governed by strongly continuous semi-groups. More precisely, for parabolic distributed systems the characterizations of regional boundary strategic sensors have been discussed and analyzed in different cases of regional boundary observability in infinite time interval. Furthermore, the results so obtained are applied in two-dimensional systems and the sensors are studied under which conditions guarantee regional boundary observability in a sub-region of the system domain boundary. Also, the authors show that, the existence of a given sensor for the diffusion system is not strategic in the usual sense, but it may be regional boundary strategic of this system.

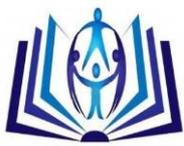

**INTRODUCTION**

The important problem of strategic sensor in distributed parameter systems has much attention in literatures ([1-2] and references therein), in order that to estimate current state of the considered system [3-4]. This problem may be called the observability notion in control systems theory [5]. Thus, the observation problem is depended on the possibility of the state reconstruction from the knowledge of system dynamics and output function by using an approach to choose the best sensor may be strategic [6-7]. Recently, regional strategic sensors characterizations is developed by El-Jai, Zerrik and Al-Saphory *et al.* for different cases in finite [8-11] or infinite time interval, may be represented as regional asymptotic systems analysis [12-16] and focused on state estimation in a sub-region $\omega$ of the domain $\Omega$ [17-18]. The purpose of this paper is to extend the previous results as in ref. [12] to the regional boundary case where the interested region $\Gamma$ is a part of the domain boundary $\partial\Omega$. The main reason behind the study of this notion is that, there exists some problem in the real world cannot observe the system state in the whole domain, but it is possible in a part of the considered domain [15-16, 18-21]. The scenario described by energy exchange problem, where the aim is to determine the energy exchanged in a casting plasma on a plane target which is perpendicular to the direction of the flow from measurements (internal pointwise sensors) carried out by thermocouples (Figure 1), where (1) is the torch of plasma, (2) is the probe of (steal), (3) is the insulator, $\Gamma$ is the face of exchange and $b_1$, $b_2$ sensor locations. This paper is organized as follows: The second section is focused on the considered system and the problem of regional boundary observability. The third section is devoted to the mathematical concepts of regional boundary observability and the characterization of regional boundary strategic sensors in various situations are studied. In the last section, we illustrate applications with many situations of sensor locations.

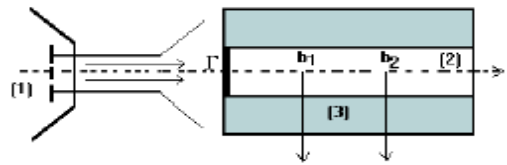

**Fig. 1:** Model of energy exchanged problem on $\Gamma$

**2. REGIONAL BOUNDARY STRATEGIC SENSORS**

In this section, we are interested to study and characterize the notion of strategic sensors on a sub-region of the domain boundary of the considered systems and present some original results related to this notion.

**2.1 Problem Statement**

Let $\Omega$ be an open regular and bounded subset of $R^n$, with smooth boundary $\partial\Omega$. Suppose that $\Gamma$ be a non-empty given sub-region of $\partial\Omega$ with positive measurement. For $T > 0$ let us set $\Theta = \Omega \times (0, \infty)$ and $\Pi = \partial\Omega \times (0, \infty)$. The considered systems is described by the following state space equations

$$\begin{cases} \frac{\partial x}{\partial t}(\xi, t) = Ax(t) + Bu(t) & \Theta \\ x(\xi, 0) = x_0(\xi) & \overline{\Omega} \\ \frac{\partial x}{\partial \vartheta}(\eta, t) = 0 & \Pi \end{cases} \quad (1)$$

where $\overline{\Omega}$ holds for closure of $\Omega$ and $x_0(\xi)$ is unknown initial state in $H^1(\overline{\Omega})$. The system (1) is defined with a Neumann boundary conditions, $\partial x/\partial \vartheta$ holds for the outward normal derivative. The measurements maybe given by the use of zone, pointwise or lines sensors which is located insides of $\Omega$ or on the boundary [1]. Thus, the augmented output function to (1) is defied by

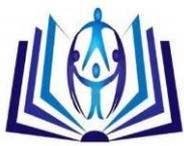

$$y(.,t) = Cx(.,t) \qquad (2)$$

where $A$ is a second order linear differential operator, which is generated a strongly continuous semi-group $(S_A(t))_{t \geq 0}$ on the Hilbert space $X = H^1(\Omega)$ and, it is self-adjoint with compact resolvent. The operator $B \in L(R^p, H^1(\Omega))$ and $C \in L(R^q, H^1(\overline{\Omega}))$, depend on the structures of actuators and sensors [1-2]. The spaces $X, U$ and $O$ be are separable Hilbert spaces where $X$ is the state space, $U = L^2(0, \infty, R^p)$ is the control space and $O = L^2(0, \infty, R^q)$ is the observation space, where $p$ and $q$ are the numbers of actuators and sensors. Under the given assumption, the system (1) has a unique solution [20]:

$$x(\xi,t) = S_A(t)x_0(\xi) + \int_0^t S_A(t-\tau)Bu(\tau)d\tau \qquad (3)$$

The problem is that, how to present sufficient conditions for regional boundary strategic sensors which enables to observe the current state in a given sub region Γ (see below Figure 2 ), using convenient sensors. Mathematical model in (Figure 2) is more general spatial case in (Figure 1).

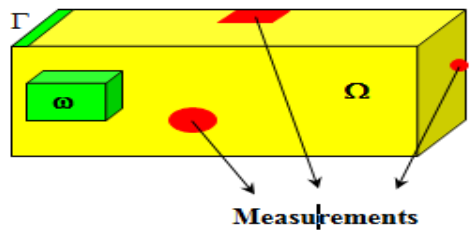

**Fig. 2:** The domain of Ω, the sub-regions ω and Γ, various sensors locations.

## 2.2 Definitions and Characterizations

The regional boundary observability concept has been developed recently by El Jai *et al.* as in [18-22] and extended to the regional boundary asymptotic state by Al-Saphory and El Jai in ref.s [1-5]. To recall regional boundary observability, consider the associated autonomous system to (1) given by

$$\begin{cases} \frac{\partial x}{\partial t}(\xi,t) = Ax(\xi,t) & \Theta \\ x(\xi,0) = x_0(\xi) & \overline{\Omega} \\ \frac{\partial x}{\partial \vartheta}(\eta,t) = 0 & \Pi \end{cases} \qquad (4)$$

Thus, the knowledge of $x(\xi,0)$ permits to observe regional boundary state $x(\xi,t)$ at any time $t$. Consider now the following points:

- The solution of system (4) is given by the following form,

$$x(\xi,t) = S_A(t)x_0(\xi), \quad \forall t \geq 0 \qquad (5)$$

- The operator $K$ is defined by following

$$\begin{aligned} K: X &\to O \\ x &\to CS_A(.)x \end{aligned} \qquad (6)$$

then, we obtain

$$y(.t) = K(t)x(.,0) \qquad (7)$$

where $K$ is bounded linear operator (this is valuable on some output function) [23].

- The operator $K^*: O \to X$ is the adjoint of $K$ defined by

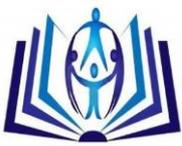

$$K^*y^* = \int_0^t S_A^*(s)\, C^*y^*(.,s)\,ds \tag{8}$$

- The trace operator of order zero

$$\gamma_0 : H^1(\Omega) \to H^{1/2}(\partial\Omega) \tag{9}$$

is linear, subjective, and continuous [3], such that $x_0^\Gamma$ is the restriction of the trace of the initial state $x_0$ to $\Gamma$. $\gamma_0^*$ denote the adjoint of $\gamma_0$ given by

$$\gamma_0^* : H^{1/2}(\partial\Omega) \to H^1(\Omega) \tag{10}$$

- For a sub-region $\Gamma \subset \partial\Omega$ and let $\chi_\Gamma$ be the restriction function defined by

$$\begin{aligned}\chi_\Gamma : H^{1/2}(\partial\Omega) &\to H^{1/2}(\Gamma) \\ x &\to \chi_\Gamma x = x|_\Gamma\end{aligned} \tag{11}$$

where $x|_\Gamma$ is the restriction of $x$ to $\Gamma$. We denote by $\chi_\Gamma^*$ the adjoint of $\chi_\Gamma$ and defined by

$$\chi_\Gamma^* : H^{1/2}(\Gamma) \to H^{1/2}(\partial\Omega) \tag{12}$$

Now, to characterize strategic sensors notion, we need some results of regional boundary observability concept in space $H^{1/2}(\Gamma)$ is extended from ref. [12, 22].

**Definition 2.1:** The system (4) augmented with the output function (2) is said to be exactly observable on $\Omega$ (or exactly $\Omega$-observable), if

$$\operatorname{Im}(K^*) = H^1(\Omega) \tag{13}$$

**Definition 2.2:** The system (4)-(2) is said to be approximately observable on $\Omega$ (or approximately $\Omega$-observable), if

$$\operatorname{Im}(\overline{K^*(.)}) = H^1(\Omega)$$

**Definition 2.3:** The system (4)-(2) is said to be regional boundary exactly observable on $\Gamma$ (or exactly $\Gamma$-observable), if

$$\operatorname{Im}(\chi_\Gamma \gamma_0 K^*) = H^{1/2}(\Gamma)$$

**Definition 2.4:** The system (4)-(2) is said to be regional boundary approximately observable on $\Gamma$ (or approximately $\Gamma$-observable), if

$$\operatorname{Im}(\overline{\chi_\Gamma \gamma_0 K^*(.)}) = H^{1/2}(\Gamma) \tag{14}$$

**Remark 2.5:** The definition 2.4 is equivalent to say that the system (4)-(2) is approximately $\Gamma$-observable if

$$\operatorname{Ker}(K(.)\gamma_0^* \chi_\Gamma^*) = \{0\} \tag{15}$$

Then, the following results can be extended from [2] to the regional boundary.

**Proposition 2.6:** The system (4)-(2) is exactly $\Gamma$-observable if there exists $\nu > 0$ such that $\forall x_0 \in H^{1/2}(\Gamma)$,

$$\|x_0\|_{H^{1/2}(\Gamma)} \leq \nu \|K\gamma_0^* \chi_\Gamma^* x_0\|_{L^2(0,\infty,\mathcal{O})} \tag{16}$$

**Proof:** The proof of this property is deduced from the usual results on observability by considering $\chi_\Gamma \gamma_0 K^*$ as in [4]. Let $E, F$ and $G$ be Banach reflexive space and $f \in L(E,G), g \in L(F,G)$, then we have

(1) $\operatorname{Im} f \subset \operatorname{Im} g$

(2) there exists $c > 0$ such that

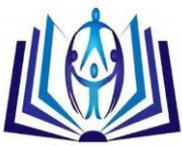

$$\|f^*x^*\|_{E^*} \le c \|g^*x^*\|_{F^*}, \quad \forall x^* \in G^*.$$

Now, if this result is applied. Choosing

$$E = G = H^{1/2}(\Gamma), \quad F = \mathcal{O}, \quad f = Id_{H^{1/2}(\Gamma)}$$

and

$$g = \chi_\Gamma \gamma_0 K^*,$$

therefore, we obtain the inequality (14) ∎.

**Corollary 2.7:** From the pervious proposition 2.6 we can get the following result:

(1) The notion of approximate Γ-observability is far less restrictive than the exact Γ-observability.

(2) From the equation (14) there exists a reconstruction error operator that gives an estimation $\tilde{x}_0$ of the initial state $x_0$ in $\Gamma$ [22]. Then, we have

$$\|x_0 - \tilde{x}_0\|_{H^{1/2}(\Gamma)} \le \|x_0 - \tilde{x}_0\|_{H^{1/2}(\partial\Omega)} \tag{17}$$

**Proposition 2.8:** The regional boundary observability concept is more convenient for the analysis of real systems [20]. Then, we can deduce that:

(1) The definitions 2.3 and 2.4 are more general and can be applied to the case where $\Gamma = \partial\Omega$.

(2) The equation (17) shows that the regional boundary state reconstruction will be more precise than if we estimate the state in the boundary of the domain $\Omega$.

(3) If a system is exactly $\Omega$-observable, then, it is exactly Γ-observable, but the converse is not true in general. Now, we prove that property (3) of remark 2.6.

**Proof:** We see that if the system is exactly observable on $\partial\Omega$, then it is exactly Γ-observable and this is a consequence of (17) and then

$$\|x_0\|_{H^{1/2}(\Gamma)} \le \|x_0\|_{H^{1/2}(\partial\Omega)}, \quad \forall x_0 \in H^{1/2}(\Gamma) \tag{18}$$

by the same way with miner tanique as in regional case [12], we can show that, if $\Gamma \subset \partial\Omega$, then

$$|x_0|_\Gamma \le |x_0|_{\partial\Omega} \tag{19}$$

and hence

$$\|x_0\|_{H^{1/2}(\Gamma)} \le \|x_0\|_{H^{1/2}(\partial\Omega)} \tag{20}$$

From equations (16), (17), (18), (19), and (20), we have

$$\| x_0 \|_{H^{1/2}(\Gamma)} = \| x_0 \|_{L(H^{1/2}(\Gamma), H^{1/2}(\partial\Omega))}$$

$$\le \|x_0\|_{H^{1/2}(\partial\Omega)}$$

$$\le \nu \|K\gamma_0^* \chi_\Gamma^* x_0\|_{L^2(0,T,\mathcal{O})}$$

Then from proposition 2.6 and remark 2.5, we can deduce that the system (4)-(2) is exactly Γ-observable.

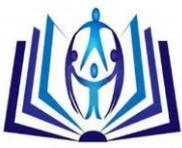

## 3. SUFFICIENT CONDITIONS FOR Γ-STRATEGIC SENSORS

The purpose of this section is to give the sufficient condition for the characterization of sensors in order that the system (1) is regionally boundary approximately observable in a region Γ.

### 3.1 Concept of Sensors

This subsection recalls and studies the concept of the sensors, which was introduced by A. El Jai [6-7]. Thus, we know that the sensors form an important link between a system and its environment [17-18]. In any case of sensors is considered via a space variable, mathematically speaking, the space variable is present in all systems described by partial differential equations [12].

**Definition 3.1:** A sensor may be defined by any couple $(D, f)$, where $D$, a non-empty closed subset of $\bar{\Omega}$, is the spatial support of sensor and $f \in L(D)$ defines the spatial distribution of the sensing measurements on $D$.

**Remark 3.2:** According to the choice of the parameters $D_i$ and $f_i$ we have various types of sensors. Sensor may be a zone types denoted by $(D_i, f_i)$, where $D_i \subset \Omega$, then, the output function (2) can be written in the form

$$y(.,t) = C\, x(.,t) = \int_{D_i} x(\xi,t) f_i(\xi) d\xi \tag{21}$$

Also, sensors maybe pointwises when $D_i = \{b_i\}$ and $f_i = \delta_{b_i}(\xi - b_i)$ represented by the couple $(b_i, \delta_{b_i})$ where $\delta_{b_i}$ is the Dirac mass concentrated in $b_i$. Thus, the output function (2) can be given by the form

$$y(.,t) = C\, x(.,t) = \int_\Omega x(\xi,t) \delta_{b_i}(\xi - b_i) d\xi \tag{22}$$

In the case of boundary zone sensors $(\Gamma_i, f_i)$ where $D_i = \Gamma_i$ with $\Gamma_i \subset \partial\Omega$ and $f_i \in L^2(\Gamma)$. Therefore, the output function (2) can then be written in the form

$$y(.,t) = C\, x(.,t) = \int_{\Gamma_i} \frac{\partial x}{\partial v}(\eta, t) f_i(\eta) d\eta \tag{23}$$

The operator $C$ is unbounded and some precautions must be taken in [3, 11].

**Definition 3.3:** A sensor $(D, f)$ is $\Omega$-strategic if the corresponding system (4)-(2) is approximately $\Omega$-observable.

**Definition 3.4:** A suit of $(D_i, f_i)_{1 \leq i \leq q}$ is said to be $\Omega$-strategic if there exists at least one sensor $(D_1, f_1)$ which is approximately $\Omega$-strategic.

**Definition 3.5:** A sensor $(D, f)$ is Γ-strategic if the corresponding system (4)-(2) is approximately Γ-observable.

**Definition 3.6:** A suit of $(D_i, f_i)_{1 \leq i \leq q}$ is said to be Γ-strategic if there exists at least one sensor $(D_1, f_1)$ which is approximately Γ-strategic.

Thus, we can deduce that the following result:

**Corollary 3.7:** A sensor is Γ-strategic if the corresponding system (4)-(2) is exactly Γ-observable.

**Proof:** Let the system (4)-(2) is exactly Γ-observable. Then, we have

$$\mathrm{Im}(\chi_\Gamma \gamma_0 K^*) = H^{1/2}(\Gamma)$$

From the decomposition sub-spaces of direct sum in Hebert space, we can represent

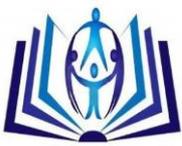

$H^{1/2}(\Omega)$ by the unique form [9]

$$Ker\ (K(t){\gamma_0}^* \chi_\Gamma^*) + Im\ (\chi_\Gamma \gamma_0 K^*) = H^{1/2}\ (\partial\Omega) \tag{24}$$

we obtain

$$Ker\ (K(t){\gamma_0}^* \chi_\Gamma^*) = \{0\}$$

This is equivalent to

$$Im\ (\overline{(\chi_\Gamma \gamma_0\ K^*(.))}) = H^{1/2}(\Gamma)$$

Finally, we can deduce this system is approximately Γ-observable and therefore this

sensor is Γ-strategic. ∎

Thus, the definition 2.3, proposition 2.6 and corollary 3.7 guarantee Γ-strategic sensors with far more restrictive conditions.

**Proposition 3.8:** From the previous results, we can deduce that:

(1) a sensor which is strategic for a system, it is Γ-strategic.

(2) a sensor which is $\Gamma_1$-strategic for a system where $\Gamma_1 \subset \partial\Omega$, then, it is $\Gamma_2$-strategic for any $\Gamma_2 \subset \Gamma_1$.

(3) One can find various sensors which are not Ω-strategic for the systems, but may be Γ-strategic and achieve the observability in Γ. This is illustrated in the following counter-example.

### 3.2 Counter-Example

Consider a two-dimensional systems described the following diffusion equations

$$\begin{cases} \frac{\partial x}{\partial t}(\xi_1, \xi_2, t) = \frac{\partial^2 x}{\partial \xi_1^2}(\xi_1, \xi_2, t) + \frac{\partial^2 x}{\partial \xi_2^2}(\xi_1, \xi_2, t) & Q \\ x(\xi_1, \xi_2, 0) = x_0(\xi_1, \xi_2) & \overline{\Omega} \\ \frac{\partial x}{\partial \vartheta}(\eta_1, \eta_2, t) = 0 & \Sigma \end{cases} \tag{25}$$

where $\Omega = [0, 1] \times [0, 1]$ and $\Gamma = [0, 1] \times \{0\}$. The output function is given by

$$y(t) = \int_{\Gamma_0} x(\eta_1, \eta_2, t)\ f(\eta_1, \eta_2)\ d\eta_1\ d\eta_2 \tag{26}$$

and $\Gamma_0 = \{0\} \times [0,1] \subset \partial\Omega$ as in (Figure 3).

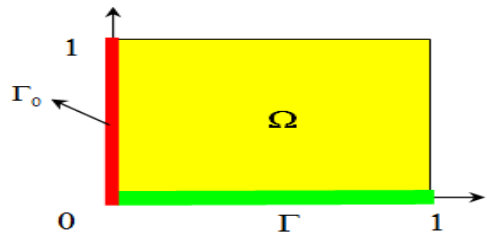

**Fig. 3 :** Domain Ω, region Γ and location $\Gamma_0$ of boundary zone sensor.

The operator $A = \frac{\partial^2}{\partial \xi_1^2} + \frac{\partial^2}{\partial \xi_2^2}$ generates a semi-group $(S(t))_{t \geq 0}$ on $H^1(\Omega)$ given by

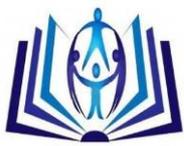

$$S(t)y = \sum_{i,j=0}^{\infty} e^{\lambda_{ij}t} \langle y, \varphi_{ij}\rangle_{H^1(\Omega)} \varphi_{ij} \qquad (27)$$

where $\lambda_{ij} = -(i^2 + j^2)\pi^2$, $\varphi_{ij}(x_1, x_2) = 2a_{ij}\cos(i\pi\xi_1)\cos(j\pi\xi_2)$ and $a_{ij} = (1 - \lambda_{ij})^{-1/2}$.

The state $x_0(\xi_1, \xi_2) = \cos(i\pi\xi_1)\cos(j\pi\xi_2)$ is not approximitely $\Omega$-observable [11]. Thus, the boundary sensor $(\Gamma_0, f)$ is not $\Omega$-strategic [1]. The systems (25)-(26) is approximately $\Gamma$-observable [22] and then boundary sensor $(\Gamma_0, f)$ is $\Gamma$-strategic [14].

In this section, we are interested to develop the results which are related to the strategic sensors and give the sufficient conditions for each sensor. For this purpose, we assume that there exists a complete set of eigenfunctions $\varphi_n$ of $A$ in $H^1(\Omega)$, associated to the eigenvalues $\lambda_n$ with a multiplicity $s_n$ and suppose that the functions $\psi_n$ defined by $\psi_n = \chi_\Gamma \gamma_0 \varphi_n$ is a complete set in $H^{1/2}(\Gamma)$ defined by, is a complete set in $H^1(\Omega)$. If the system (2.1) has $J$ unstable modes, then we have the following result.

**Theorem 3.7:** Assume that $sup\ r_n = r < \infty$, then the suite $(\Gamma_i, f_i)_{1 \leq i \leq q}$ of zones boundary sensors are $\Gamma$-strategic if and only if

(1) $q \geq r$

(2) $rank\ G_n = r_n$, where

$$G_n = (G_n)_{ij}\ \text{with}\ 1 \leq i \leq q, 1 \leq j \leq r_n\ \text{and}\ (G_n)\ \text{is given by}$$

$$(G_n)_{ij} = \begin{bmatrix} \langle \varphi_{n_1}, f_1(.)\rangle_{L^2_{(\Gamma_1)}}, & \ldots, & \langle \varphi_{n_{r_n}}, f_1(.)\rangle_{L^2_{(\Gamma_1)}} \\ & \vdots & \\ \langle \varphi_{n_1 1}, f_q(.)\rangle_{L^2_{(\Gamma_q)}}, & \ldots, & \langle \varphi_{n_{r_n}}, f_q(.)\rangle_{L^2_{(\Gamma_q)}} \end{bmatrix}$$

**Proof:** The proof is developed in the case where the suit of sensors are of boundary zones type $(\Gamma_i, f_i)_{1 \leq i \leq q}$ and located on $\partial \Omega$. If the suit of sensors are $\Gamma$-strategic, then the corresponding system (25)-(26) is approximately $\Gamma$-observable [7], it is equivalent to

$$[K\gamma_0^* \chi_\Gamma^* x^* = 0 \Longrightarrow x^* = 0],\ \text{for}\ x^* \in H^{1/2}(\Gamma) \qquad (28)$$

We have

$$K\gamma_0^* \chi_\Gamma^* x^* = (\sum_{n\geq 1} \exp^{(\lambda_n t)} \sum_{j=1}^{r_n} \langle \varphi_{nj}, \gamma_0^* \chi_\Gamma^* x^*\rangle_{H^1(\Omega)} \langle \varphi_{nj}, f_i\rangle_{L^2(\Gamma)})_{1\leq i \leq q}$$
$$= (\sum_{n\geq 1} \exp^{(\lambda_n t)} \sum_{j=1}^{r_n} \langle \chi_\Gamma \gamma_0 \varphi_{nj}, x^*\rangle_{H^{1/2}(\Gamma)} \langle \varphi_{nj}, f_i\rangle_{L^2(\Gamma)})_{1\leq i \leq q}$$
$$= (\sum_{n\geq 1} \exp^{(\lambda_n t)} \sum_{j=1}^{r_n} \langle \psi_{nj}, x^*\rangle_{H^{1/2}(\Gamma)} \langle \varphi_{nj}, f_i\rangle_{L^2(\Gamma)})_{1\leq i \leq q}$$

Now, if the suit of sensors is not strategic sensors, then the system (25)-(26) is not approximately $\Gamma$-observable, and hence there exists $x^* \neq 0$ such that

$$K\gamma_0^* x_\Gamma^* x^* = 0 \Longleftrightarrow \sum_{j=1}^{r_n} \langle \psi_{nj}, x^*\rangle_{H^{1/2}(\Gamma)} \langle \varphi_{nj}, f_i\rangle_{L^2(\Gamma)} = 0\ \forall n, n \geq 1$$

Suppose that $x_n$ defined by

$$x_n = \begin{bmatrix} \langle \psi_{n_1}, x^*\rangle_{H^{1/2}(\Gamma)} \\ \vdots \\ \langle \psi_{n_{r_n}}, x^*\rangle_{H^{1/2}(\Gamma)} \end{bmatrix} \qquad (29)$$

Then

$$G_n x_n = 0,\ \forall n \geq 1 \Longleftrightarrow rank\ G_n \neq r_n.$$

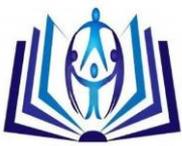

Conversely, if rank $G_n \neq r_n$ for some $n$, then there exists

$$x_n = \begin{bmatrix} x_{n_1} \\ \vdots \\ x_{n_{r_n}} \end{bmatrix} \neq 0, \quad x^* = \sum_{j=1}^{r_n} x_{nj} \psi_{nj} \in H^{1/2}(\Gamma) \neq 0$$

such that

$$G_n x_n = 0$$

Thus, we can find a non-zero $x^* \in H^{1/2}(\Gamma)$ such that

$$\langle x^*, \psi_{jk} \rangle_{H^{1/2}(\Gamma)} = 0 \text{ if } j \neq n$$

and

$$\langle x^*, \psi_{nk} \rangle_{H^{1/2}(\Gamma)} = x_{nk}, 1 \leq k \leq r_n$$

For which $\sum_{k=1}^{r_j} \langle \varphi_{jk}, f_i \rangle_{\Gamma_i} \langle x^*, \psi_{jk} \rangle_{H^{1/2}(\Gamma)} = 0 \quad j \neq n, \ 1 \leq i \leq q$

and also

$$\sum_{k=1}^{r_j} \langle \varphi_{nk}, f_i \rangle_{\Gamma_i} \langle x^*, \psi_{nk} \rangle_{H^{1/2}(\Gamma)} = 0, 1 \leq i \leq q$$

otherwise there exists $x^* \neq 0 \in H^{1/2}(\Gamma)$, such that

$$K \gamma_0^* x_\Gamma^* x^* = 0,$$

Thus, the system (25)-(26) is not approximately Γ-observable and then the sensors are not Γ-strategic. ∎

**Corollary 3.8:** If the system (25)-(26) is exactly Γ-observable, then, the rank condition in theorem 3.7 is satisfied.

**Remark 3.9:** The previous result can be extended to the case of internal zone, filament and internal or boundary sensors as in ref.s [12-16].

**Remark 3.10:** The important to introduce this notion is that the using to charaterize the regional boundary exponential reduced observability in distributed parameter system as in [24] and this notion is extended to mutipule situations for finite time interval [25-26] or infinite as in [27-29].

### 4. APPLICATION TO SENSOR LOCATIONS

In this section, we present an application of the above results in two-dimensional systems defined on $\Omega = [0, a_1] \times [0, a_2]$ by the form

$$\begin{cases} \frac{\partial x}{\partial t}(\xi_1, \xi_2, t) = \Delta x(\xi_1, \xi_2, t) & \Theta \\ x(\xi_1, \xi_2, 0) = x_0(\xi_1, \xi_2) & \overline{\Omega} \\ \frac{\partial x}{\partial \vartheta}(\eta_1, \eta_2, 0) = 0 & \Pi \end{cases} \quad (30)$$

together with output function is described by (2). Let $\Gamma = \{a_1\} \times (0, a_2)$ be the considered region is subset of $[0, a_1] \times [0, a_2]$. In this case, the eigenfunctions of system (30) are given by

$$\varphi_{ij}(\xi_1, \xi_2) = \frac{2}{\sqrt{a_1 a_2}} \cos i\pi \left(\frac{\xi_1}{a_1}\right) \cos j\pi \left(\frac{\xi_2}{a_2}\right) \quad (31)$$

associated with eigenvalues

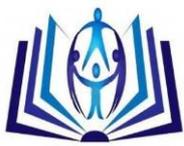

$$\lambda_{ij} = -\left(\frac{i^2}{(a_1)^2} + \frac{j^2}{(a_2)^2}\right) \quad (32)$$

The following results give information on the locations of internal, boundary zone or pointwise Γ-strategic sensors.

### 4.1 Zone Sensor Cases

This subsection study various types of domains with different systems.

#### 4.1.1 Rectangular domain

We discuss and examine different type of zone sensors.

***Internal rectangular zone case:***

Consider the system (30)-(2) where the sensor supports $D$ are located inside Ω. Then the output (2) can be written by the form

$$y(t) = \int_D x(\xi_1, \xi_2, t) f_i(\xi_1, \xi_2) d\xi_1 d\xi_2 \quad (33)$$

where $D \subset \Omega$ is location of zone sensor and $f \in L^2(D)$. In this case of (Figure 4), the eigenfunctions and the eigenvalues

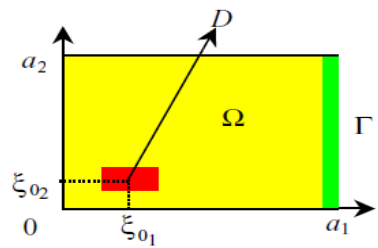

**Fig. 4:** Domain Ω, region Γ and location $D$ of internal zone sensor.

are given by (31) and (32). However, if we suppose that

$$\frac{(a_1)^2}{(a_2)^2} \notin \mathbb{N} \quad (34)$$

where $\mathbb{N}$ is the natural numbers. If $r = 1$ then one sensor $(D, f)$ maybe suffices to achieve Γ-strategic sensor of the corresponding systems (30)-(33) [9-11]. Let the measurement support is rectangular with

$$D = [\xi_1 - l, \xi_1 + l_1] \times [\xi_2 - l_2, \xi_2 + l_2] \in \Omega$$

then, we have the following result.

**Corollary 4.1:** If $f_1$ is symmetric about $\xi_1 = \xi_{01}$ and $f_2$ is symmetric about $\xi_2 = \xi_{02}$, then the sensor $(D, f)$ is Γ-strategic to the systems (30)-(33) if $i(\xi_{01})/(a_1)$ and $j(\xi_{02})/(a_2) \notin \mathbb{N}$ for every $i, j = 1, \ldots, J$

***One side boundary zone case:***

In the case where $\Gamma_0 \subset \partial\Omega$ and $f \in L^2(\Gamma_0)$, the sensor $(\Gamma_0, f)$ *is* located on the one side of the boundary $\partial\Omega$ in $\Gamma_0 = [(\eta_{01} - l, (\eta_{02} - l] \times \{a_2\}$ as in (Figure 5). Consider again the systems (30)-(33), then the output function is given by

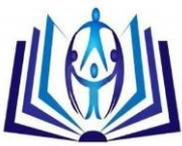

$$y(t) = \int_{\Gamma_0} x(\eta_1,\eta_2,t)\, f(\eta_1,\eta_2)\, d\eta_1\, d\eta_2 \qquad (35)$$

then, we obtain.

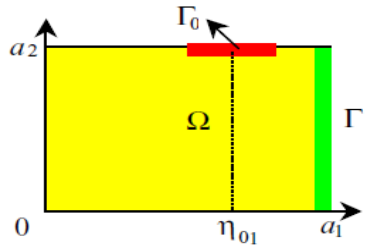

**Fig. 5:** Domain $\Omega$, region $\Gamma$ and location $\Gamma_0$ of boundary on one side zone sensor.

**Corollary 4.2:** Suppose that $f$ is symmetric with respect to $\eta_1 = \eta_{01}$, then the sensor $(\Gamma_0, f)$ is $\Gamma$-strategic to the systems (30)-(35) if $i\,(\eta_{01})/(a_1) \notin \mathbb{N}$ for every $i = 1, \dots, J$.

***Two sides boundary zone case:***

Suppose that the sensor $(\bar\Gamma, f)$ is located on $\bar\Gamma = [0, \bar\eta_{01} + l_1] \times \{0\} \cup \{0\} \times [0, \bar\eta_{02} + l_2] = \bar\Gamma_1 \cup \bar\Gamma_2 \subset \partial\Omega$ and $f_{\bar\Gamma_1}$ is symmetric with respect to $\bar\eta_1 = \bar\eta_{01}$ and the function $f_{\bar\Gamma_2}$ is symmetric with respect to $\bar\eta_2 = \bar\eta_{02}$ as in (Figure 6). Thus, the output function is given by

$$y(t) = \int_{\bar\Gamma} x(\eta_1,\eta_2,t)\, f(\eta_1,\eta_2)\, d\eta_1\, d\eta_2 \qquad (36)$$

then we have.

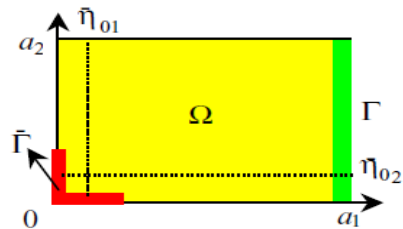

**Fig. 6:** Domain $\Omega$, region $\Gamma$ and location $\bar\Gamma$ of two sides zone sensor.

**Corollary 4.3:** If $f_1$ is symmetric about $\bar\eta_1 = \bar\eta_{01}$ and $f_2$ is symmetric about $\bar\eta_2 = \bar\eta_{02}$ then the sensor $(\Gamma, f)$ is $\Gamma$-strategic to the systems (30)-(36) if $i\,(\bar\eta_{01})/(a_1)$ and $j\,(\bar\eta_{02})/(a_2) \notin \mathbb{N}$ for every $i,j = 1, \dots, J$.

### 4.1.2 Disc domain

We explore some results concern different type of zone sensors in disc domain.

***Internal circular zone case:***

In this case, systems (30) may be given by the following form

$$\begin{cases} \frac{\partial x}{\partial t}(r,\theta,t) = \Delta x(\xi_1,\xi_2,t) & \Theta \\ x(r,\theta,0) = x_0(r,\theta) & \bar\Omega \\ \frac{\partial x}{\partial \vartheta}(a,\theta,t) = 0 & \Pi \end{cases} \qquad (37)$$

where $0 < \theta < 2\pi$, $\Omega = (0,a)$, $r = a > 0$, and $\theta \in [0, 2\pi]$, $t > 0$ are defined as in (Figure 7). The augmented output function described by

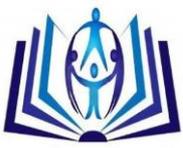

$$y_i(t) = \int_{D_i} x(r,\theta,t) f_i(r_i,\theta_i) dr_i d\theta_i \tag{38}$$

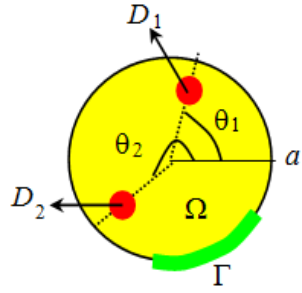

**Fig. 7:** Disc domain $\Omega$, region $\Gamma$ and locations $D_1, D_2$ of internal pointwise sensors.

Let the eigenfunctions and eigenvalues concerning the region $\Gamma = (a, \theta_i)_{2 \leq i \leq q}$ of $\partial\Omega$ with $\in [0, 2\pi]$ are defined

$$\lambda_{ij} = \beta_{ij}^2, i \geq 0, j \geq 1 \tag{39}$$

where $\beta_{ij}$ $\beta_{nm}$ are the zeros of the Bessel functions $J_n$ and

$$\begin{cases} \psi_{ij}(r,\theta) = J_0(\beta_{ij}^2 r), & j \geq 1 \\ \psi_{ij}(r,\theta) = J_n(\beta_{ij_1}^2 r) \cos(i\theta) & i, j_1 \geq 1 \\ \psi_{ij}(r,\theta) = J_n(\beta_{ij_2}^2 r) \sin(i\theta\} & i, j_2 \geq 1 \end{cases} \tag{40}$$

with multiplicity $S_{ij} = 2$ for all $ij \neq 0$ and $S_{ij} = 1$ for all $ij = 0$. In this case, the $\Gamma$-strategic sensor is required at least two zone sensors $(D_i, f_i)_{2 \leq i \leq q}$ where $D_i = (r_i, \theta_i), i = 1,2$ (see [14]). If we consider the case of Dirichlet or mixed boundary conditions, we can get various functions [2]. Thus, we develop some practical examples by using the symmetry conditions. If $f_i$ and $D_i$ are symmetric with respect to $\theta = \theta_i$, for all $2 \leq i \leq q$, then we have .

**Corollary 4.4:** If the sensors $(D_i, f_i)_{2 \leq i \leq q}$ are located in $D_i = (r_i, \theta_i), i = 1,2$ and $i(\theta_1 - \theta_2)/\pi \notin \mathbb{N}$ for every $i, j = 1, \ldots, J$, then $(D_1, f_1)$ and $(D_2, f_2)$ are is $\Gamma$-strategic to the systems (37)-(38).

**Boundary circular zone case:**

In this case the system (30) is augmented with output function described by

$$y_i(t) = \int_{\partial\Omega} x(r,\theta,t) \delta_{\Gamma_i}(r - r_i, \theta - \theta_i) dr_i d\theta_i \tag{41}$$

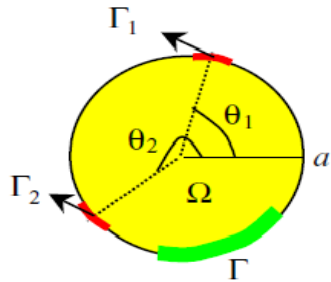

**Fig. 8:** Disc domain $\Omega$, region $\Gamma$ and locations $\Gamma_1, \Gamma_2$ of internal pointwise sensors.

When the sensors $(\Gamma_i, f)_{2 \leq i \leq q}$ are located on $\partial\Omega$ and the function $f_{\Gamma_i}$ is symmetric with respect to $\theta = \theta_i$, for all $2 \leq i \leq q$, as in (Figure 8). So, we have.

**Corollary 4.5:** If the sensors $(\Gamma_i, f_i)_{2 \leq i \leq q}$ are located in $\Gamma_i = (a, \theta_i), i = 1, 2$ and $i(\theta_1 - \theta_2)/\pi \notin \mathbb{N}$ for every $i, j = 1, \dots, J$, then $(\Gamma_1, f_1)$ and $(\Gamma_2, f_2)$ are $\Gamma$-strategic to the systems (37)-(41).

### 4.2 Pointwise Sensor Cases

This subsection study again various types of domains in different systems.

#### 4.2.1 The domain $\Omega = [0, a_1] \times [0, a_2]$

We deal with different type of pointwise sensors.

*Internal pointwise case:*

Let us consider the case of pointwise sensor located inside of $\Omega$. Thus, the system (30) is augmented with the following output function.

$$y(t) = \int_\Omega x(\xi_1, \xi_2, t)\delta(\xi_1 - b_1, \xi_2 - b_2)d\xi_1 d\xi_2 \qquad (42)$$

where $b = (b_1, b_2)$ is the location of pointwise sensor as defined in (Figure 9)

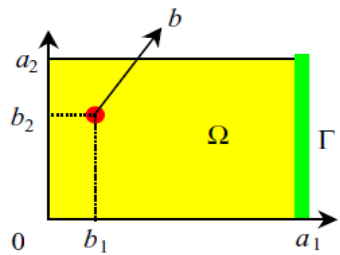

**Fig. 9:** The domain $\Omega$, region $\Gamma$ and location $b$ of internal pointwise sensor.

In this case may be one pointwise sensor $(b, \delta_b)$ is sufficient for strategic sensor on $\Gamma$ of systems (30)-(42). Thus, we obtain the following result.

**Corollary 4.6:** The sensor $(b, \delta_b)$ is $\Gamma$-strategic to the systems (30)-(42) if $i(b_1)/(a_1)$ and $j(b_2)/(a_2) \notin \mathbb{N}$, for every $i, j = 1, \dots, J$.

*Internal filament case:*

Consider the case where the information is given on the curve $\sigma = Im(\gamma)$ with $\gamma \in C^1(0, 1)$ (Figure 10).

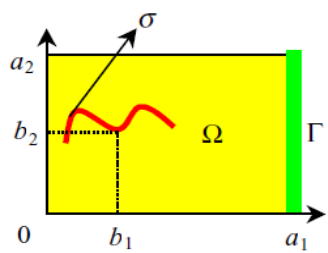

**Fig. 10:** The domain $\Omega$, and location $\sigma$ of internal filament sensors.

If the measurements recovered by filament sensor $(\sigma, \sigma_b)$ such that is symmetric with respect to the line $b = b_i$, then we have.

**Corollary 4.7:** The sensor $(\sigma, \sigma_b)$ is $\Gamma$-strategic to the systems (30)-(42) if $i(b_1)/i(a_1)$ and $i(b_2)/i(a_2) \notin \mathbb{N}$, for every $i$.

*Boundary pointwise case:*

Now, the system (30) is augmented with the following output function.

$$y(t) = \int_{\partial\Omega} x(\xi_1, \xi_2, t)\delta\,(0, \xi_2 - b_2)d\xi_2 \qquad (43)$$

where $b = (0, b_2)$ is the location of pointwise sensor $(b, \delta_b)$ as defined in (Figure 11).

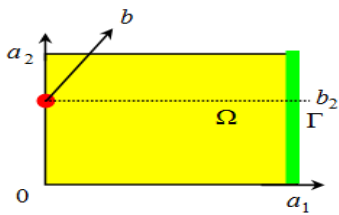

**Fig. 11:** The domain $\Omega$, and location $b$ of boundary pointwise sensor.

Thus, we obtain the following result.

**Corollary 4.8:** The sensor $(b, \delta_b)$ is $\Gamma$-strategic to the systems (30)-(43) if $i(b_2)/(a_2) \notin \mathbb{N}$, for every $i = 1, \ldots, J$.

**4.2.2 The domain $\Omega = [0, a_1] \times [0, a_2]$**

We discuss and examine different type of zone sensors.

*Internal pointwise case:*

Here, the system (37) is augmented with the following output function.

$$y_i(t) = \int_{\Omega} x(r, \theta, t)\delta_{p_i}(r - r_i, \theta - \theta_i)dr_i d\theta_i \qquad (44)$$

where $0 < \theta < 2\pi$, $\Omega = (0, a)$, and $r = a > 0$. The locations of pointwise sensors $(p_i, \delta_{p_i})_{2 \le i \le q}$ are $(p_1, \delta_{p_1})$, $(p_2, \delta_{p_2})$ with $p_1 = (r_1, \theta_1)$ and $p_2 = (r_2, \theta_2)$ in $\Omega$ (Figure 12), then we can get the following result.

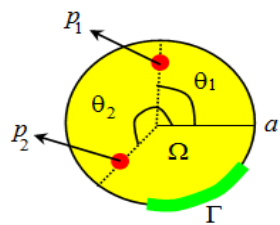

**Fig. 12:** Disc domain $\Omega$, region $\Gamma$ and locations $p_1, p_2$ of internal pointwise sensors.

**Corollary 4.9:** The sensor $(p_i, \delta_{p_i})_{2 \le i \le q}$ are located in $p_i = (r_i, \theta_i) \in \Omega$ is $\Gamma$-strategic to the systems (37)-(44) if $i(\theta_1 - \theta_2)/\pi \notin \mathbb{N}$ for every $i, j = 1, \ldots, J$

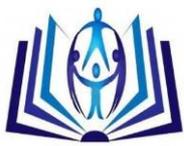

*Boundary pointwise case:*

Here, the system (37) is augmented with the following output function.

$$y_i(t) = \int_{\partial\Omega} x(r,\theta,t)\delta_{p_i}(r - r_i, \theta - \theta_i) dr_i d\theta_i \qquad (45)$$

where $0 < \theta < 2\pi$, $\Omega = (0,a)$, and $r = a > 0$. The locations of pointwise sensors $(p_1, \delta_{p_1})$, $(p_2, \delta_{p_2})$ with $p_1 = (a, \theta_1)$ and $p_2 = (a, \theta_2)$ in $\partial\Omega$ (Figure 13), then we can get.

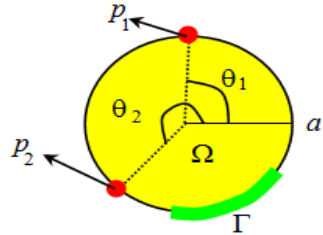

**Fig. 13:** Disc domain $\Omega$, region $\Gamma$ and locations $p_1, p_2$ of boundary pointwise sensors.

**Corollary 4.10:** The sensors $(p_i, \delta_{p_i})_{2 \leq i \leq q}$ are located in $p_i = (a, \theta_i) \in \partial\Omega$ are $\Gamma$-strategic to the systems (37)-(45) if $i(\theta_1 - \theta_2)/\pi \notin \mathbb{N}$ for every $i, j = 1, \dots, J$.

**Corollary 4.11:** These results can be extended to the following:

1. Case of Dirichlet or mixed boundary conditions [1-6].

2. We can show that the observation error decreases when the number and support of sensors increases [11].

## 5. CONCLUSION

The notion of regional boundary strategic sensors have been developed and examined. A various regional boundary observability have been discussed and analyzed which permit us to avoid some bad sensor locations. In addition, many interesting results concerning the choice of such sensors are given and illustrated in specific situations with diffusion systems. Thus, several questions still opened, for example, the simulations of this model are under consideration and the problem of finding an optimal sensor location ensuring such an objective.

**ACKNOWLEDGMENTS** Our thanks in advance to the editors and experts for considering this paper to publish in this estimated journal and for their efforts.

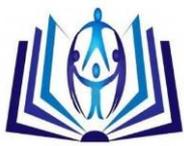

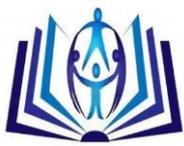

**Author's Biographies**:

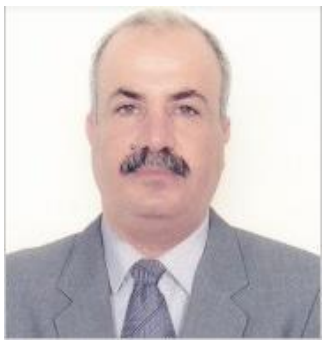

Raheam Al-Saphory is full professor at the TIKRIT University in IRAQ. He received his Ph. D. degree in Control System and Analysis in (2001) from LTS/ Perpignan University France. He has a post doctorate as a researcher in 2001-2002 at LTS. Al-Saphory wrote one book and many papers in the area of systems analysis and control. Also he is a supervisor of many Ph D. and M.Sc. students and he was Ex-head of Department of Mathematics /College of Education for Pure Science Tikrit University 2010-2011. He visits many Centers and Scientific Departments of Bangor University/ Wales/ UK with academic staff of Iraqi Universities in 2013. Now, he is a chermam of Scientific and postgraduate studies committee at the Department of Mathematics 2014-present.

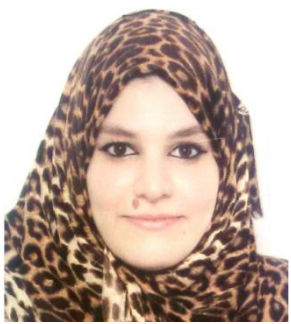

Hind K. Kolaib is a researcher and she received here Ms. c.. degree in Control System and Analysis in (2017) from the Department of Mathematics/ College of Education for Pure Sciences/ TIKRIT University / IRAQ. Here research area focused on Distributed Parameter Systems Analysis and Control.